\documentclass[12pt]{article}
\usepackage{amssymb}
\usepackage{amsmath}
\usepackage[usenames]{color}
\usepackage{graphicx}
\usepackage{amsthm}

\usepackage{amsmath,amssymb,mathtools}
\usepackage{amsmath,amssymb,mathtools,amsthm}
	\usepackage{textcomp}
	\usepackage{amsfonts}
	\usepackage{graphicx}
	\usepackage{enumerate}
\usepackage{color}
\usepackage{times}
\usepackage[utf8]{inputenc}
\usepackage[T1]{fontenc}
\usepackage{mathtools}
\usepackage{geometry}
\usepackage{bm}
\newcommand{\ncom}{\newcommand}
\newcommand{\non}{\nonumber}
\newcommand{\noi}{\noindent}
\ncom{\vone}{\vskip 2ex}
\ncom{\vtwo}{\vskip 4ex}

\newtheorem{remark}{Remark}
\ncom{\vsa}{\vspace{.3cm}}
\ncom{\vsb}{\vspace{.4cm}}

\begin{document}

\title{Generalization and Probabilistic Proofs of Some Combinatorial Identities }

%\date{} for date
\author{Palaniappan Vellaisamy and Puja Pandey\\
\small Department of Statistics and Applied Probability \\[-0.8ex]
\small University of California Santa Barbara \\[-0.8ex] 
\small Santa Barbara, CA, 93106, USA. \\
\small \text{Email}:   pvellais@ucsb.edu; pujapandey@ucsb.edu}
\maketitle

\begin{abstract}
  Using a probabilistic approach, we derive some interesting combinatorial identities involving gamma and beta functions. These results generalize certain well-known combinatorial identities involving binomial coefficients and special functions. In particular, by studying moments of the difference of two gamma and beta random variables, both in the dependent and independent cases, we obtain new combinatorial identities. This approach provides a systematic method to derive further combinatorial identities from probabilistic transformations.

\end{abstract}
\section{Introduction}
There are several interesting combinatorial identities involving binomial coefficients, gamma functions and hypergeometric functions  (see, for example, Riordan \cite{Riordan},  Bagdasaryan \cite{Bagdasaryan},  Vellaisamy \cite{Vellaisamy}, and {the references therein). One of these is the famous identity that involves the convolution of the central binomial coefficients: \begin{equation} \label{eqn1.1}
 \sum_{k=1}^{n} \binom{2k}{k}\binom{2n-2k}{n-k}=4^n.
\end{equation}   
In recent years, researchers have provided several proofs of \eqref{eqn1.1}. 
A proof that uses generating functions can be found in Stanley \cite{Stanley}. Combinatorial proofs can also be found, for example, in Sved \cite{Sved}, De Angelis \cite{De Angelis}, and Miki{\' c} \cite{Mickic}. 
\noindent A related and an interesting identity for the alternating convolution of central binomial coefficients is 
\begin{equation}\label{e12}
	\sum_{k=0}^{n} (-1)^k \binom{2k}{k}  \binom{2n-2k}{n-k} =
	\begin{cases*}
		2^n \binom{n}{\frac{n}{2}}, & if  $n$ is even  \\
		\phantom{0 }0, & if  $n$ is odd.
	\end{cases*} 
\end{equation}

\noindent The combinatorial proofs of the above identity can be found in Nagy \cite{Nagy}, Spivey \cite{Spivey}, and  Miki{\' c} \cite{Mickic}.
Recently, there has been considerable interest in finding simple probabilistic proofs for combinatorial identities; see, for example, Vellaisamy and Zeleke \cite{VellaisamyZeleke} and the references therein. A probabilistic proof of the  identity in \eqref{e12} can be found in Pathak \cite{Pathak}. Chang and  Xu \cite{ChungXu} extended the result in \eqref{eqn1.1} and presented a probabilistic proof of the identity
\begin{equation}\label{eqn13}
 \sum_{\substack{k_j \geq 0,\; 1 \leq j \leq m;\\ \sum_{J=1}^m  k_j=n}} \binom{2k_1}{k_1}\binom{2k_2}{k_2}\cdots\binom{2k_m}{k_m}=\frac{4^n}{n!} \frac{\Gamma(n+\frac{m}{2})}{\Gamma(\frac{m}{2})},
\end{equation}
 for positive integers $m$ and $n.$ Miki{\' c} \cite{Mickic} discussed a combinatorial proof of \eqref{eqn13} based on the method of recurrence relations and telescoping.
 Vellaisamy and Zeleke \cite{VellZel19} provided a  probabilistic connection between Euler's constant and the Basel problem. Also, Vellaisamy and Zeleke \cite{VellZel19a} presented a probabilistic proofs several identities involving
 beta functions.
 
 \noindent Another interesting combinatorial identity  (see (6.60) in Gould (2010))  is
 \begin{equation}\label{e14}
 \frac{1}{4^{2n}}\sum_{k=0}^{2n} (-1)^k \binom{2n}{k} \binom{2k}{k} \binom{4n-2k}{2n-k}= \frac{1}{4^{2n}} \binom{2n}{n}^2. 
\end{equation}
 
  Our goal in this paper is  to derive some combinatorial identities similar to the ones given in
   in \eqref{e12}-\eqref{e14}, using a simple probabilistic approach.  Indeed, we derive certain identities involving beta and gamma functions. Our method uses the difference of two gamma or beta random variables,
   and their moments. Computing moments in two different ways lead to 
   some interesting combinatorial identities, which generalize some
   known identities. One of the main goals of this paper is to provide
   a probabilistic proof of the identity in \eqref{e14}.

\section{ Identities Derived from  Gamma Distributions}
\noindent Let $ X_1$ and $X_2$ be two independent standard gamma  random variables with shape parameter $p >0$ and scale parameter unity,  denoted by Ga$(p)$. Note that the density of a random variable   $X \sim \text{Ga}(p)$ is 
$$f(x|p)=\frac{1}{\Gamma(p)}e^{- x}x^{p-1},~~x>0,\; p >0.$$
 It is known that (see Rohatgi and Saleh \cite{RohatgiSaleh}), for $n \geq 1$, 
\begin{equation*}
E(X^n)=\frac{\Gamma(n+p)}{\Gamma(p)},
\end{equation*} 
where $\Gamma(p)$ denotes the usual gamma function.

\noindent Also, it is well known that $(X_1+X_2)$ has a gamma distribution with parameter $2p$, that is, $( X_1+X_2)  \sim \text{Ga}(2p)$ and so
\begin{equation} \label{eqn25}
E(X_1+X_2)^{2n}=\frac{\Gamma(2p+2n)}{\Gamma(2p)}.
\end{equation}
 Next we compute the even moments of $( X_1-X_2)$. Note that the exact distribution of $( X_1-X_2)$ is too complicated to use to find the
 moments (see Mathai \cite{Mathai}). So, we use the moment generating function approach.\\

\noindent It is known (see Rohatgi and Saleh \cite{RohatgiSaleh}) that the  $MGF$ of $X_1$ is $M_{X_1}(t)= E(e^{tX_1})= (1-t)^{-p}.$ Hence, the $MGF$ of $X_1-X_2$ is

\vspace*{-.8cm}
\begin{align}
M_{X_1-X_2}(t) =& M_{X_1}(t)M_{X_2}(-t) \nonumber \\
=& (1-t)^{-p} (1+t)^{-p} 
= (1-t^2)^{-p},        
\end{align}
which exists for $|t|<1.$ 

\noindent Also, for $p>0$ {and}  $|q|<1$, it is known that
\begin{equation*}
(1-q)^{-p}= \sum_{n=0}^{\infty} \frac{\Gamma(n+p) q^n}{\Gamma(n+1) \Gamma(p)}.
\end{equation*}
Using the above result, we have
\begin{equation*} 
M_{X_1-X_2}(t)= (1-t^2)^{-p}= \sum_{n=0}^{\infty} \frac{\Gamma(n+p) t^{2n}}{\Gamma(n+1) \Gamma(p)}.
\end{equation*}
Hence, for  $n \geq 1$, 
\begin{equation} \label{eqn27}
E(X_1-X_2)^{2n}=M_{X}^{(2n)}(t)|_{t=0}=  \frac{\Gamma(n+p) {(2n)!}}{\Gamma(n+1) \Gamma(p)},
\end{equation}	
where $f^{(k)}$ denotes the $k$-th derivative of $f$.

\noindent Therefore,  from \eqref{eqn25} and \eqref{eqn27}, we get 
\begin{equation} \label{eqn28}
\frac{1}{2} E\Big((X_1+X_2)^{2n}+(X_1-X_2)^{2n}\Big)=\frac{1}{2} \Bigg(\frac{\Gamma(2p+2n)}{\Gamma(2p)}+\frac{\Gamma(p+n)(2n)!}{\Gamma(p)\Gamma(n+1)}\Bigg).
\end{equation}
\noindent On the other hand, using the binomial theorem and $ (-1)^{2n-k}= (-1)^k$, we get 
\begin{align} \label{eqn29}
\frac{1}{2} E\Big((X_1+X_2)^{2n}+(X_1-X_2)^{2n}\Big)=& \frac{1}{2} E \Bigg(\sum_{k=0}^{2n} \binom{2n}{k} X_1^k  X_2^{2n-k}\nonumber\\ 
+&  \sum_{k=0}^{2n} (-1)^k\binom{2n}{k} X_1^k  X_2^{2n-k} \Bigg) \nonumber \\
=&\frac{1}{2}E \Bigg(2 \sum_{k=0}^{n} \binom{2n}{2k} X_1^{2k}  X_2^{2n-2k} \Bigg)\nonumber\\
=& \sum_{k=0}^{n} \binom{2n}{2k} E(X_1)^{2k} E (X_2)^{2n-2k}\nonumber \\
=& \sum_{k=0}^{n} \binom{2n}{2k}\frac{\Gamma(2k+p)}{\Gamma(p)}\frac{\Gamma(2n-2k+p)}{\Gamma(p)}.
\end{align}
for all $p>0$ and  $ n \geq 1. $
\noindent From \eqref{eqn28} and \eqref{eqn29}, we get an identity involving the gamma functions 
\begin{align} \label{eqn6}
\sum_{k=0}^{n} \binom{2n}{2k}\frac{\Gamma(2k+p)}{\Gamma(p)}\frac{\Gamma(2n-2k+p)}{\Gamma(p)}=&
\frac{1}{2}\Bigg(\frac{\Gamma(2p+2n)}{\Gamma(2p)}+ \frac{\Gamma(n+p)(2n)!}{\Gamma(n+1)\Gamma(p)}\Bigg)
\end{align}
for all $p>0$ and $ n \geq 1.$
\subsection{The Special Case $p=\frac{1}{2}$}
The identity 
\begin{equation*} 
\sum_{k=0}^{n} \binom{4k}{2k}\binom{4n-4k}{2n-2k}=2^{4n-1} +  2^{2n-1}\binom{2n}{n}
\end{equation*}
first appeared in Brychkov \cite{Brychkov}. A proof of this identity can be found in Vignat and Moll \cite{VigMoll15}. Here we will show that this identity follows as a special case from the identity in \eqref{eqn6}. 
Also, we use  the identity 
\begin{align} \label{eqn11}
\frac{\Gamma(n+\frac{1}{2})}{\Gamma(\frac{1}{2})} = \frac{\binom{2n}{n} n!}{4^n}, 
\end{align}
which can easily be verified.

\noindent Let now $p=\frac{1}{2}.$ Then
the right hand side (rhs) of equation \eqref{eqn6} reduces to, 
\begin{align} \label{eqn7}
 \frac{1}{2}\Bigg(\frac{\Gamma(2n+1)}{\Gamma(1)}
+\frac{\Gamma(n+\frac{1}{2})}{\Gamma(n+1)}\frac{(2n)!)}{\Gamma(\frac{1}{2})}\Bigg)\nonumber\\
=&\frac{1}{2}\Bigg((2n)!+\frac{\Gamma(n+\frac{1}{2})}{\Gamma(\frac{1}{2}) n!}(2n)!\Bigg)\nonumber \\
=&\frac{1}{2}\Bigg((2n)!+\frac{\binom{2n}{n} n!}{4^n} \frac{(2n)!}{n!}\Bigg)\nonumber \\
=&\frac{1}{2}\Bigg((2n)!+\frac{\binom{2n}{n}}{4^n} {(2n)!}\Bigg)
\end{align}
Next, the left hand side (lhs) of equation \eqref{eqn6}, when  $p=\frac{1}{2},$ 
becomes 
\begin{align} \label{eqn8}
\sum_{k=0}^{2n} \binom{2n}{2k}\frac{\Gamma(2k+\frac{1}{2})}{\Gamma(\frac{1}{2})}\frac{\Gamma(2n-2k+\frac{1}{2})}{\Gamma(\frac{1}{2})} &  \nonumber \\
=&\sum_{k=0}^{n} \binom{2n}{2k}\frac{\binom{4k}{2k} (2k)!}{4^{2k}} \binom{4n-4k}{2n-2k} \frac{(2n-2k)!}{4^{2n-2k}} \nonumber \\
=&
\sum_{k=0}^{n} \binom{2n}{2k}\binom{4k}{2k}\binom{4n-4k}{2n-2k} \frac{ (2k)! (2n-2k)!}{4^{2n}}\nonumber \\
=& (2n)! \sum_{k=0}^{n} \binom{4k}{2k}\binom{4n-4k}{2n-2k}\frac{1}{4^{2n}}
\end{align}
Equating \eqref{eqn7} and \eqref{eqn8}, we get, 
$$(2n)! \sum_{k=0}^{n} \binom{4k}{2k}\binom{4n-4k}{2n-2k}\frac{1}{4^{2n}}=\frac{(2n)!}{2}\;\Bigg[1+\frac{\binom{2n}{n}}{4^n}\; \Bigg] \nonumber.$$
That is, 
\begin{align} \label{eqn9}
\sum_{k=0}^{n} \binom{4k}{2k}\binom{4n-4k}{2n-2k}=& \frac{4^{2n}}{2} \;\Bigg[1+\frac{\binom{2n}{n}}{4^n}\; \Bigg] \nonumber \\
=& \frac{4^{n}}{2} \;\Bigg[2^{2n}+\binom{2n}{n}\; \Bigg] \nonumber \\
=& 2^{4n-1} +  2^{2n-1}\binom{2n}{n}, 
\end{align}
which is the identity (7.5) in Gould \cite{Gould}.
See also the identity  (3.1) of Vignat and Moll \cite{VigMoll15}.

\begin{remark} {\em
To prove the above identity, Vignat and Moll \cite{VigMoll15} computed the
even moments $(X_1-X_2)$, where $X_1$ and $X_2$ follow $Ga(\frac{1}{2})$
random variables, using a different approach. They used the facts that
\begin{align}
(X_1-X_2)\stackrel{d}{=} (X_1+X_2) Y,
\end{align}
where $Y$ is independent of $X_1$ and $X_2$ and has the symmetric beta density  given by
\begin{align}
f(y)= \frac{1}{\pi \sqrt{1-y^2}}, ~ |x|<1
\end{align}
and
\begin{align}
E(Y^{2n})= \frac{1}{2^{2n}} \binom{2n}{n}.
\end{align}
Note that we have computed the even moments of  $(X_1-X_2),$ 
where $X_1$ and $X_2$ follow $Ga(p)$ variables, using the direct
moment generating function approach.} 
\end{remark}

\subsection{Special Case $p=n$}
Consider now the special case  $p=n.$  Then the  identity in \eqref{eqn6} 
becomes
\begin{align*}
\sum_{k=0}^{n} \binom{2n}{2k} \frac{\Gamma(n+2k)}{\Gamma (n)}\frac{\Gamma(3n-2k)}{\Gamma (n)}=&\frac{1}{2} \Bigg( \frac{\Gamma(4n)}{\Gamma(2n)}
+\frac{\Gamma(2n)\Gamma(2n+1)}{\Gamma(n)\Gamma(n+1)}\Bigg)
\end{align*}
which leads to an interesting identity
\begin{align*}
\sum_{k=0}^{n} \binom{2n}{2k} \Gamma(n+2k){\Gamma(3n-2k)}=&\frac{\Gamma^2(n)}{2} \Bigg( \frac{\Gamma(4n)}{\Gamma(2n)}
+2\frac{\Gamma^2(2n)}{\Gamma^2(n)}\Bigg).
\end{align*}

\subsection{Another Identity involving Beta  Functions}
\noindent Consider next the moments of $T= \frac{X_1-X_2}{X_1+X_2}$, where 
$X_1$ and $X_2$ are independent Ga$(p)$ random variables. Note first that
\begin{align*}
 E(T^{2n})= E( 1-2Y_2)^{2n},
\end{align*}
where $Y_2= \frac{X_2}{X_1+X_2} \sim \text{Be}(p, p)$, the beta  distribution, as is well known.\\

\noindent Using binomial theorem
and the fact $ (-1)^{2n-k}= (-1)^{k}$, we get
\begin{align}\label{eqn18}
E( 1-2Y_2)^{2n} =& \frac{1}{B(p, p)} \int_{0}^{1} (2y-1)^{2n} y^{p-1} (1-y)^{p-1} dy \nonumber \\
    =& \frac{1}{B(p, p)} \sum_{k=0}^{2n} (-1)^k \binom{2n}{k}
     2^k \int_{0}^{1}  y^{p+k-1} (1-y)^{p-1} dy \nonumber\\
=& \frac{1}{B(p, p)} \sum_{k=0}^{2n} (-1)^k \binom{2n}{k}
2^k B(p+k, p)\nonumber \\
=&  \sum_{k=0}^{2n} (-1)^k \binom{2n}{k}
2^k     \frac{B(p+k, p)}{B(p, p)},
\end{align}
where $B(p, p) $ denotes the usual beta function.

\vspace{0.3cm}
\noindent Next we compute the $2n$-th moment of $T$ using its density. 
 Note that the density $f(t)$ of $T$ is symmetric about origin and so
suffices to find it for some $ 0 <t < 1.$ Using the fact
\begin{align*}
P(T \leq -t)=& P( 1-2Y_2 \leq -t)= P\Big(Y_2 \geq  \frac{1+t}{2} \Big)\\
\implies F(-t)    =& \frac{1}{B(p,p)} \int_{\frac{1+t}{2}}^{1} y^{p-1} (1-y)^{p-1}.
\end{align*}
Differentiating with respect to  $t$, we get
the density of $T$ as
\begin{align*}
f_T(t) =& \frac{1}{B(p, p) 2^{2p-1}} (1-t^2)^{p-1} \\ 
=& \frac{1}{B(\frac{1}{2}, p)} (1-t^2)^{p-1},~~ -1<t<1,
\end{align*} 
since $ B(p, p)2^{2p-1} = B(\frac{1}{2}, p)$, which follows 
from, for $p>0$,
\begin{align*}
\frac{\Gamma{(p+ \frac{1}{2})}}{\Gamma{(\frac{1}{2})}}=\frac{\Gamma{(2p+1)}} {\Gamma{(p+1)}} \frac{1}{4^p} \cdot
\end{align*} 

Hence,
\begin{align} \label{eqn19}
E(T^{2n})=& \frac{1}{B(\frac{1}{2}, p)}  \int_{-1}^{1}  t^{2n} (1-t^2)^{p-1} dt \nonumber \\
=& \frac{2}{B(\frac{1}{2}, p)}  \int_{0}^{1}  t^{2n} (1-t^2)^{p-1} dt \nonumber \\
=& \frac{1}{B(\frac{1}{2}, p)}  \int_{0}^{1}  y^{n} (1-y)^{p-1} y^{-\frac{1}{2}}dy \nonumber \\
=& \frac{1}{B(\frac{1}{2}, p)}  \int_{0}^{1}  y^{n+\frac{1}{2}-1} (1-y)^{p-1} dy \nonumber \\
=& \frac{B(n+\frac{1}{2}, p)}{B(\frac{1}{2}, p)}
\end{align}

\noindent From \eqref{eqn18} and \eqref{eqn19}, we get another interesting   identity 
\begin{align} \label{eqn20}
 \sum_{k=0}^{2n} (-2)^k \binom{2n}{k}
     \frac{B(p+k, p)}{B(p, p)} =& \frac{B(n+\frac{1}{2}, p)}{B(\frac{1}{2}. p)},
\end{align}
for all $p>0$ and $ n \geq 1.$

\noindent Consider next the special case $p= \frac{1}{2}.$ The   lhs of \eqref{eqn20}
becomes
\begin{align} \label{eqn21}
\sum_{k=0}^{2n} (-2)^k \binom{2n}{k}
\frac{B(p+k, p)}{B(p, p)}=&\sum_{k=0}^{2n} (-2)^k \binom{2n}{k}
\frac{\Gamma (k+\frac{1}{2})}{\Gamma (\frac{1}{2})}\frac{1}{\Gamma{(k+1)}} \nonumber \\
 =& \sum_{k=0}^{2n} (-2)^k \binom{2n}{k}
 \frac {\binom{2k}{k}}{4^k} k!  \frac{1}{\Gamma{(k+1)}}   \nonumber \\
 =& \sum_{k=0}^{2n} (-1)^k \binom{2n}{k}
\binom{2k}{k}\frac{1}{2^k} \cdot 
\end{align}
Note that in the second step above, we have used the identity in \eqref{eqn11} .

\noi Similarly, the  rhs of \eqref{eqn20} becomes, when $p= \frac{1}{2},$ becomes
\begin{align} \label{eqn22}
\frac{B(n+\frac{1}{2}, p)}{B(\frac{1}{2}, \frac{1}{2})}=&  \frac {\Gamma{(n+\frac{1}{2})}} {\Gamma{(\frac{1}{2})}} \frac{1}{\Gamma{(n+1)}}  \nonumber \\
=&\frac {\binom{2n}{n}}{4^n} n!  \frac{1}{\Gamma{(n+1)}}   \nonumber \\
=& \frac {\binom{2n}{n}}{4^n}.
\end{align}

\noindent Thus, from  \eqref{eqn21}  and  \eqref{eqn22}, we get an interesting combinatorial identity
\begin{align} \label{eqn23}
\sum_{k=0}^{2n} (-1)^k \binom{2n}{k}
\binom{2k}{k}\frac{1}{2^k}=& \frac {\binom{2n}{n}}{4^n}.
\end{align}
or to put it another way,

\begin{align} \label{eqn24}
\sum_{k=0}^{2n} (-1)^k \binom{2n}{k}
\binom{2k}{k}{2^{2n-k}}=& \binom{2n}{n}
\end{align}
which is a new combinatorial identity. \\

\noi It is interesting to note that the derivation of the above  combinatorial identity  relies on the $2n$-th moment of $T$, where $T = \frac{X_1}{X_1 + X_2} - \frac{X_2}{X_1 + X_2}$.
Clearly, $\frac{X_1}{X_1 + X_2}$ and $\frac{X_2}{X_1 + X_2}$ are dependent random variables, each following a Be$(p,p)$ distribution. Thus, the identity is obtained by considering the difference of two dependent beta random variables. What happens when we consider the moments of two  independent beta random variables? We address this question in the next section.

\section{A Probabilistic Proof of the Identity in (4)}
Let $Y_1$ and $Y_2$ be now independent $\text{Be}(\frac{1}{2},\frac{1}{2})$ variables and  
$X=Y_1-Y_2.$ Then, using \eqref{eqn11},
\begin{align} \label{e27} 
E(X)^{2n}=& \sum_{k=0}^{2n} (-1)^k \binom{2n}{k} E(Y_{1}^k)E(Y_{2}^{2n-k}) \nonumber \\
=&\sum_{k=0}^{2n} (-1)^k \binom{2n}{k} \frac{B(k+\frac{1}{2}, \frac{1}{2})}{B(\frac{1}{2}, \frac{1}{2})}  \frac{B(2n-k+\frac{1}{2}, \frac{1}{2})}{B(\frac{1}{2}, \frac{1}{2})}  \nonumber \\
=& \sum_{k=0}^{2n} (-1)^k \binom{2n}{k} \frac{\Gamma(k+\frac{1}{2})}{\Gamma(\frac{1}{2})} \frac{1}{k!} \frac{\Gamma(2n-k+\frac{1}{2})}{\Gamma(\frac{1}{2})}  \frac{1}{(2n-k)!}  \nonumber\\
=& \sum_{k=0}^{2n} (-1)^k \binom{2n}{k} \frac{\binom{2k}{k}} {4^k}
    \binom{4n-2k}{2n-k} \frac{1}{4^{2n-k}}  \nonumber \\
=& \frac{1}{4^{2n}}\sum_{k=0}^{2n} (-1)^k \binom{2n}{k} \binom{2k}{k} \binom{4n-2k}{2n-k}.     
\end{align}

\noi Next, we compute $E(X)^{2n}$ using the density of $X$. But, a direct computation of rather tedious computations lead to an expression which is similar or reduces  to the rhs of \eqref{e27}. So, it was  quite a challenging one, as can be seen in the following derivations, to 
get the expression in the rhs of \eqref{e14}.  

\vone
\noi Let $F_{1}$ denote the Appell’s first hypergeometric function of two variables defined by
$$F_{1}(a,b_{1}, b_{2};c;x_{1},x_{2}) = \sum_{k=0}^{\infty} \sum_{m= 0}^{\infty} \frac{(a_{1},m+k)}{(c,m+k)} (b_{1},m) (b_{2},k) \frac{x_{1}^{m} x_{2}^{k}}{m! \; k!},$$
\normalsize
where $(a,m)$ denotes the Pochhammer symbol, defined by
$$(a,m) = a (a+1)...(a+m-1) = \frac{\Gamma(a+m)}{\Gamma(a)}, m >0.$$

\noi When $Y_{1}$ and $Y_{2}$ follow the 
$\text{Be}(\frac{1}{2},\frac{1}{2})$ distribution, the density function of $X = Y_{1} - Y_{2}$ is given by (see \cite{PGT93})
\[
f_{X}(x) = \begin{cases}
	\frac{1}{B(\frac{1}{2}, \frac{1}{2})} F_{1}(\frac{1}{2},0, \frac{1}{2}; 1; 1-x, 1-x^{2}), & \text{for }~ 0 \leq x <1, \\
	\frac{1}{B(\frac{1}{2}, \frac{1}{2})} F_{1}(\frac{1}{2}, \frac{1}{2},0; 1; 1-x^{2}, 1+x), & \text{for}~ -1 \leq x < 0.
\end{cases}
\]
Note that, using the series expression of $F_1$,
\begin{align*}
	F_{1}(\frac{1}{2},0, \frac{1}{2}; 1; 1-x, 1-x^{2}) = \sum_{k=0}^{\infty} \sum_{m= 0}^{\infty} \frac{(\frac{1}{2},m+k)}{(1,m+k)} (0,m) (\frac{1}{2},k) \frac{(1-x)^{m} (1-x^{2})^{k}}{m! \; k!},
\end{align*}
where the sum over $m$ vanishes for $m >0$ since $(0,m)=0, m >0$ and $(0,0) = 1$. Thus,
\begin{align*}
	F_{1}(\frac{1}{2},0, \frac{1}{2}; 1; 1-x, 1-x^{2}) &= \sum_{k=0}^{\infty} \frac{(\frac{1}{2},k)  (\frac{1}{2},k)}{(1,k)}  \frac{ (1-x^{2})^{k}}{k!} \\
	&= {}_2F_1(\frac{1}{2}, \frac{1}{2};1;1-x^{2}),
\end{align*}
where ${}_2F_1(a,b;c;z) = \sum_{k=0}^{\infty} \frac{(a,k)(b,k)}{(c,k)} \frac{z^{k}}{k!}$ is the Gauss-hypergeometric function. \\
Similarly, $F_{1}(\frac{1}{2}, \frac{1}{2},0; 1; 1-x^{2}, 1+x)$ also reduces to  ${}_2F_1(\frac{1}{2}, \frac{1}{2};1;1-x^{2})$.
Thus, $X$ have a symmetric density over $(-1, 1).$ This implies, $E(X^{2n+1})=0$ for all $n \ge 1.$ Next, we compute the even moments of $X$.
For $n \ge 1$,
\begin{align}\label{e28}
	\mathbb{E}[X^{2n}] &= \frac{2}{\pi} \int_{0}^{1} x^{2n} {}_2F_1(\frac{1}{2}, \frac{1}{2};1;1-x^{2})~dx. 
\end{align}
Making the substitution $x^{2} = u$ in \eqref{e28}, we get
\begin{align}\label{e29}
	\mathbb{E}[X^{2n}] &= \frac{1}{\pi} \int_{0}^{1} u^{n-\frac{1}{2}} {}_2F_1(\frac{1}{2}, \frac{1}{2};1;1-u)~du. 
\end{align}
Applying the Euler-type integral representation of the hypergeometric series ${}_2F_1(a,b;c;z)$:
$$B(b,c-b) \; {}_2F_1(a,b;c;z) = \int_{0}^{1} y^{b-1} (1-y)^{c-b-1} (1-zy)^{-a} ~dy,$$
we get
\begin{align*} 
	\mathbb{E}[X^{2n}] &= \frac{1}{\pi^{2}} \int_{0}^{1} u^{n-\frac{1}{2}} \int_{0}^{1} y^{-\frac{1}{2}} (1-y)^{-\frac{1}{2}} (1-(1-u)y)^{-\frac{1}{2}}~dy \; du \nonumber \\
	&= \frac{1}{\pi^{2}} \int_{0}^{1} \int_{0}^{1} u^{n-\frac{1}{2}}  y^{-\frac{1}{2}} (1-y)^{-\frac{1}{2}} (1-(1-u)y)^{-\frac{1}2{}}~dy \; du
\end{align*} 
Making the substitution $ v= 1- (1-u)y$ in the above integral gives
\begin{align*}
	\mathbb{E}[X^{2n}] &= \frac{1}{\pi^{2}} \int_{u =0}^{1} \int_{v=u}^{1} u^{n-\frac{1}{2}} \; v^{-\frac{1}{2}} \; (1-v)^{-\frac{1}{2}} \; (v-u)^{-\frac{1}{2}} ~dv \; du.
\end{align*}
Now interchanging the order of integration gives
\begin{align*}
	\mathbb{E}[X^{2n}] &= \frac{1}{\pi^{2}} \int_{v=0}^{1} \int_{u=0}^{v} u^{n-\frac{1}{2}} \; v^{-\frac{1}{2}} \; (1-v)^{-\frac{1}{2}} \; (v-u)^{-\frac{1}{2}}~du \; dv.
\end{align*}
Upon applying the substitution $u=vs$ in the above equation, the integral takes the form
\begin{align}\label{e30}
	\mathbb{E}[X^{2n}] &= \frac{1}{\pi^{2}} \int_{v=0}^{1} \int_{s=0}^{1} (vs)^{n-\frac{1}{2}} \; v^{-\frac{1}{2}} \; (1-v)^{-\frac{1}{2}} \; (v-vs)^{-\frac{1}{2}}v~ds \; dv \non \\
	&= \frac{1}{\pi^{2}} \int_{v=0}^{1} \int_{s=0}^{1} v^{n-\frac{1}{2}} \; s^{n-\frac{1}{2}} \; (1-v)^{-\frac{1}{2}} \; (1-s)^{-\frac{1}{2}}  \; ds \; dv  \non\\
	&= \frac{1}{\pi^{2}} \int_{v=0}^{1} v^{n-\frac{1}{2}} \; (1-v)^{-\frac{1}{2}} \; dv \int_{s=0}^{1}  s^{n-\frac{1}{2}} \;  (1-s)^{-\frac{1}{2}}  \; ds  \non\\
	&=  \frac{1}{\pi^{2}} \Big( B\left(n+\frac{1}{2},\frac{1}{2}\right)\Big)^{2}  \non\\
    &= \frac{1}{4^{2n}} \binom{2n}{n}^{2},
\end{align}
using \eqref{eqn22} and $B(\frac{1}{2}, \frac{1}{2})= \pi.$

\noi From \eqref{e27} and \eqref{e30}, the identity in \eqref{e14} is probabilistically proved.

\vone
\noi \textbf{Remark 2.} We have shown
\begin{equation}\label{e31}
\mathbb{E}[X^{2n}]=	\frac{1}{\pi} \int_{0}^{1} u^{n-\frac{1}{2}} \; {}_2F_1(\frac{1}{2}, \frac{1}{2};1;1-u) \;du =  \frac{1}{4^{2n}} \binom{2n}{n}^{2}.
\end{equation}
Using the series representation of the hypergeometric function, the integral on the lhs of \eqref{e31} can be further expressed as
\begin{align}\label{e32}
	\frac{1}{\pi} \int_{0}^{1} u^{n-\frac{1}{2}} \; {}_2F_1(\frac{1}{2}, \frac{1}{2};1;1-u) \;du  &=  \frac{1}{\pi} \int_{0}^{1} u^{n-\frac{1}{2}} \; \sum_{k=0}^{\infty}  \frac{(\frac{1}{2},k)^{2}}{(1,k)} (1-u)^{k}\;du  \non \\
	&= \frac{1}{\pi} \sum_{k=0}^{\infty}  \frac{(\frac{1}{2},k)^{2}}{(1,k)} \int_{0}^{1} u^{n-\frac{1}{2}} (1-u)^{k}\;du \non\\
	&= \frac{1}{\pi} \sum_{k=0}^{\infty}  \frac{(\frac{1}{2},k)^{2}}{k!} B\left(n+\frac{1}{2},k+1\right) \non\\
	&= \frac{1}{\pi} \sum_{k=0}^{\infty}  \frac{(\frac{1}{2},k)^{2} \Gamma(n+\frac{1}{2})}{\Gamma(n+k+\frac{3}{2})} 
\end{align}
Consequently, the final reduced form of \eqref{e32} provides the following identity:
$$\sum_{k=0}^{\infty}  \frac{(\frac{1}{2},k)^{2} \; \Gamma(n+\frac{1}{2})}{\Gamma(n+k+\frac{3}{2})} =   \frac{\pi}{4^{2n}} \binom{2n}{n}^{2}.$$

\end{document}